\documentclass[reqno]{amsart}
	\usepackage{cmbright}
\usepackage[hidelinks]{hyperref}
\usepackage{doi,url}
\usepackage{enumitem}
\usepackage{tikz}

\usepackage[utf8]{inputenc}
\usepackage{amsmath}
\usepackage{amssymb} 
\usepackage[numbers]{natbib}
\usepackage{prettyref}
\newtheorem{theorem}{Theorem}

\theoremstyle{plain}
\makeatletter
\newtheorem*{thm@P}{Theorem \pipp@}

\makeatother

\newenvironment{pf}{\begin{proof}}{\end{proof}}

\numberwithin{equation}{section}

\newtheorem*{thm*}{Theorem}
\newrefformat{thm}{Theorem~\ref{#1}} 
\makeatletter 
\let\old@newtheorem\newtheorem 
\renewcommand{\newtheorem}[2]{\old@newtheorem{#1}[thm]{#2} 
	\newrefformat{#1}{#2~\ref{##1}}} 
\makeatother 

\theoremstyle{definition}

\theoremstyle{remark}

\newcommand{\N}{\mathbb N}

\newcommand{\PA}{\text{PA}}
\newcommand{\ZF}{\text{ZF}}
\newcommand{\Con}{\text{Con}}

\title[]{Is the twin prime conjecture independent of Peano Arithmetic?}
\date{1 Nov. 2021}
\author[A. Berarducci]{Alessandro Berarducci}
\address{Università di Pisa, Dipartimento di Matematica, Largo Bruno Pontecorvo 5, 56127 Pisa}
\author[A. Fornasiero]{Antongiulio Fornasiero}
\address{Università di Firenze,
	Dipartimento di Matematica e Informatica "Ulisse Dini",
	viale Morgagni 67/a, 50134 Firenze}
\author[J. D. Hamkins]{Joel David Hamkins}
\address{University of Oxford, University College, High Street, Oxford OX1 4BH}
\urladdr{http://jdh.hamkins.org}

\begin{document}
\begin{abstract}
	We show that there is an arithmetical formula $\varphi$ such that $\ZF$ proves that $\varphi$ is independent of $\PA$ and yet, unlike other arithmetical independent statements, the truth value of $\varphi$ cannot at present be established in $\ZF$ or in any other trusted metatheory. In fact we can choose an example of such a formula $\varphi$ such that $\ZF$ proves that $\varphi$ is equivalent to the twin prime conjecture. We conclude with a discussion of notion of trustworthy theory and a sharper version of the result. 
\end{abstract}
\maketitle

We are looking for an arithmetical formula $\varphi$ whose truth value (in the standard model $\N$) is at present unknown and such that $\ZF$ (Zermelo-Fraenkel set theory) proves that $\varphi$ is independent of $\PA$ (first order Peano Arithmetic). 

A solution of the riddle has been obtained independently by the first two authors in a recent arXiv preprint and by the third author in an older post on Mathoverflow, which actually provides a sharper version given in Theorem \ref{thm:rosser}

We take the opportunity to reflect on the nature of the incompleteness results and make some philosophical considerations.

Let us first observe that the formula $\Con(\ZF)$ expressing the consistency
of $\ZF$ would not work. Although $\Con(\ZF)$ is independent
of~$\PA$, the proof of this fact cannot be carried out in $\ZF$. To see why,
recall that an independence result involves both the unprovability of the
formula and the unprovability of its negation. In our case one direction poses
no problems: $\ZF$ does prove that $\Con (\ZF)$ is unprovable in $\PA$. The
problem lies in the other direction: $\ZF$, if consistent, does not prove that
$\lnot \Con (\ZF)$ is unprovable in $\PA$. This depends on the fact that
$\lnot \Con (\ZF)$ has complexity $\Sigma^0_1$, so if true it must be provable
in~$\PA$, and if unprovable it must be false, provably in $\ZF$ (and also in
$\PA$ itself). Thus if $\ZF$ proves that $\lnot \Con(\ZF)$ is unprovable
in~$\PA$, then $\ZF$ proves $\Con(\ZF)$, contradicting G\"odel's second
incompleteness theorem.

More generally, we cannot take for $\varphi$ any formula of complexity
$\Pi^0_1$ (such as Goldbach's conjecture). The argument is the same: if $\ZF$
proves that a $\Pi^0_1$ formula $\phi$ is independent of $\PA$, then $\ZF$ proves
$\phi$, so we would have established the truth value of $\phi$ on the basis of
a trusted metatheory, which is against our desiderata. By simmetry any formula
of complexity $\Sigma^0_1$ does not work.

Before presenting the solution to our riddle, the reader is invited to make an
attempt himself. In the meantime we try to clarify the notion of trustworthy
metatheory. Let us agree that if we trust a theory $T$ we must at least
believe that $T$ does not prove false arithmetical sentences, or, said
differently, we must believe that the standard model $\N$ of $\PA$ is also a
model of the arithmetical consequences of $T$. Now, if we trust $\ZF$, we must
believe that $\ZF$ is consistent (as otherwise $\ZF$ would prove $0=1$, a
false arithmetical statement), hence we must trust $\ZF+\Con(\ZF)$, because
$\Con(\ZF)$ is true in $\N$ if and only if $\ZF$ is consistent. We can now
iterate the process: if we trust $\ZF+\Con(\ZF)$ we must also trust
$\ZF+ \Con(\ZF+\Con(\ZF))$ and so on.  This line of thought is at the origin
of the notion of transfinite progression of theories explored by Turing
\cite{Turing1939} and Feferman \cite{Feferman1962}. Although an infinite
iteration of this process is somewhat delicate, it must be admitted that the
finite iterations, starting from $\ZF$, are generally considered
trustworthy. There is little doubt that a proof of Riemann's hypothesis in
$\ZF+\Con(\ZF)$ would count as a genuine proof of Riemann's hypothesis.

Returning to our problem, recall that we want an arithmetical formula $\varphi$, such that $\ZF$ proves that $\varphi$ is independent of $\PA$ and such that truth value of $\varphi$ is at present unkonwn in $\ZF$ or in any other trustworthy metatheory. This is a challenge, since it is difficult to prove anything about a formula $\varphi$, even if only its independence, without having a clue as to whether it is true or false. Nevertheless we can argue as follows. Working inside the metatheory $\ZF$, take a true arithmetical sentence $\alpha$ independent of $\PA$, for instance $\alpha = \Con(\PA)$. Let $\beta$ be another true arithmetical sentence independent of $\PA+\alpha$, for instance $\beta = \Con(\PA+\alpha)$. Now let $\gamma$ be an arbitrary arithmetical statement, for instance a $\Pi^0_2$-formula expressing the existence of infinitely many twin primes. Our final formula $\varphi$ is the implication $\alpha \to (\beta \land \gamma)$. Note that $\varphi$ is true if and only if $\gamma$ is true (since $\alpha$ and $\beta$ are true). Let us show that $\varphi$ is independent of $\PA$.
If $\PA$ proves $\varphi$, then in particular $PA+\alpha$ proves $\beta$, contradicting the independence of $\beta$ from $PA+\alpha$. If on the other hand $\PA$ proves $\lnot \varphi$, then,  tautologically, $\PA$ proves $\alpha$, contradicting the independence of $\alpha$ from $\PA$.

We have thus solved the riddle: in $\ZF$ the formula $\varphi$ is independent of $\PA$ and yet it is equivalent to the twin prime conjecture, whose truth value is at present unknown in any trustworthy theory. 

One may be tempted to conclude that the twin prime conjecture is independent
of $\PA$, a formidable result, but this would be cheating, because the formula
$\varphi$ is an unnatural way of expressing the twin prime conjecture in $\PA$
($\varphi$ is provably equivalent to the ``natural formalization'' of the twin prime
conjecture in $\ZF$, not in $\PA$). This leads to the question of what is a
natural formalization and which are the means to identify~it.  In
\cite{Feferman1960} Feferman taught us that intensionally correct translations
must preserve logical complexity. For instance for the validity of G\"odel's
second incompleteness theorem it is important that the axioms of $\PA$ are
presented by a formula of complexity $\Sigma^0_1$, as otherwise there are
counterexamples (e.g.\ one may describe the axioms of $\PA$ by the
$\Sigma^0_2$ formula representing the union of all the consistent fragments
$\text{I}\Sigma^0_n$ of $\PA$). However in our case complexity alone cannot be the
answer, since our formula $\varphi$ has complexity $\Pi^0_2$, exactly as one would
expect from a natural formalization of the twin prime conjecture. Having
discarded complexity, we have no proposals other than to resort to the notion
of intended meaning, however vague it may be.

It is worth to observe that our formula $\varphi$ would remain undecidable in $\PA$, provably in $\ZF$, even if somebody found a proof of the twin prime conjecture in $\PA$. In this event however the truth value of $\varphi$ would be decided in a trustworthy metatheory, specifically in $\ZF$, so $\varphi$ would no longer satisfy our desiderata.
This raises the problem of whether there could be undecidable arithmetical formulas unaffected by such unforeseeable future event. The correct formulation of the question would however require further thoughts. The notion of absolutely undecidable mathematical statement has been investigated in \cite{Koelnner2010}, but the candidates there considered are not arithmetical.


Before presenting the sharper version of the result mentioned at the beginning of the paper, let us show that the informal notion of trustworthy theory is somewhat problematic. We have already remarked that a necessary condition is that a trustworthy theory is arithmetically sound, namely it does not prove false arithmetical statements. Being primarily interested in arithmetical statements, we may be tempted to take arithmetical soundness as the formal counterpart, formalizable in $\ZF$, of the notion of trustworthy theory. This is however at odds with our intuition. To see why, consider the theory $\ZF^\# = \ZF + \text{``}\ZF \text{ is not arithmetically sound}\text{''}$, where the statement between quotes can be formalized as a $\Delta^1_1$-formula. Nobody in his right mind would consider $\ZF^\#$ as a trustworthy theory (one would rather trust $\ZF+ \text{``$\ZF$ is arithmetically sound''}$) and yet it is not difficult to prove, in $\ZF$, that $\ZF^\#$ is arithmetically sound if so is $\ZF$.
Indeed, reason in the metatheory $\ZF+ \text{``$\ZF$ is arithmetically sound''}$ and suppose for a contradiction that there is a true arithmetical formula $\psi$ such that $\ZF^\#\vdash \lnot \psi^\omega$, where $\psi^\omega$ is the translations of $\psi$ in the language of $\ZF$.
Then $\ZF+ \psi^\omega \vdash \text{``$ZF$ is arithmetically sound''}$. This can be strengthened to $\ZF+ \psi^\omega \vdash \text{``$ZF+\psi^\omega$ is arithmetically sound''}$ (because $\ZF$ proves that an arithmetically sound theory remains such with the addition of a true arithmetical formula). In particular $\ZF+\psi^\omega \vdash \Con(\ZF + \psi^\omega)$. Thus, by G\"odel's second incompleteness theorem, $\ZF \vdash \lnot \psi^\omega$, contradicting the assumption that $\ZF$ is arithmetically sound and $\psi$ is a true arithmetical formula.

Informal notions like ``natural formalization" or ``trustworthy theory'', can serve as guiding principles but can hardly be formalized without betraying their meaning.

A use of the Rosser sentence enables the following slightly sharper version of the result, proved by the third author on MathOverflow \cite{Hamkins}, where the assumption that the theory is trustworty is weakened to mere consistency. The claim is that for every statement $\psi$, such as the twin primes conjecture, there is a variant way to express it, $\varphi$, which is equivalent to the original statement $\psi$, but which is formally independent of any particular desired consistent theory $T$ interpreting $\PA$. 

\begin{theorem}\label{thm:rosser}
Suppose that $\psi$ is any arithmetical sentence, such as the twin primes conjecture, and
$T$ is any consistent recursively axiomatized theory interpreting $\PA$. 
Then there is another arithmetical sentence  $\varphi$ such that
\begin{enumerate}
\item $\text{PA}+\text{Con}(T)$ proves that $\psi$ and $\varphi$ are equivalent.
\item $T$ does not prove $\varphi$.
\item $T$ does not prove $\neg\varphi$.
\end{enumerate}
\end{theorem}

\begin{pf}
	Let $\rho$ be the Rosser sentence for $T$, the self-referential assertion that for any proof of $\rho$ in $T$, there is a smaller proof of $\neg\rho$. The Gödel-Rosser theorem establishes that if $T$ is consistent, then $T$ proves neither $\rho$ nor $\neg\rho$. Formalizing this shows that $\text{PA}+\text{Con}(T)$ proves that $\rho$ is not provable in $T$ and hence that $\rho$ is (vacuously) true and that $\text{PA}$ proves $\text{Con}(T)\to\text{Con}(T+\rho)$. By the incompleteness theorem applied to $T+\rho$, therefore, $T+\rho$ does not prove $\text{Con}(T)$. In short, $T+\rho$ is a strictly intermediate theory between $T$ and $T+\text{Con}(T)$.
	
	Now, fix any statement $\psi$ and let $\varphi$ be the assertion $\rho\to (\text{Con}(T)\wedge \psi)$. Since $\text{PA}+\text{Con}(T)$ proves $\rho$, elementary logic shows that $\text{PA}+\text{Con}(T)$ proves that $\psi$ and $\varphi$ are equivalent. The statement $\varphi$, however, is not provable in $T$, since if it were, then $T+\rho$ would prove $\text{Con}(T)$, which it does not by our observations above. Conversely, $\varphi$ is not refutable in $T$, since any such refutation would mean that $T$ proves that the hypothesis of $\varphi$ is true and the conclusion false; in particular, it would require $T$ to prove the Rosser sentence $\rho$, which it does not by the Gödel-Rosser theorem.
\end{pf}\goodbreak

Thus, what we require of trustworthiness in our theory is merely that $\Con(T)$ is true. And of course any instance of non-provability from a theory $T$ will require the consistency of $T$, and so one cannot provide a solution to the problem without assuming the theory is consistent.

\bibliographystyle{is-alpha}

\end{document}